\documentclass[10pt]{article}
\usepackage{amssymb}
\usepackage{amsfonts}
\title{The small index property for countable superatomic boolean algebras}
\author{J K Truss} 
\date{University of Leeds}%\footnotemark}

\begin{document}
\maketitle  
\newtheorem{lemma}{Lemma}[section]
\newtheorem{theorem}[lemma]{Theorem}
\newtheorem{corollary}[lemma]{Corollary}

\setcounter{footnote}{1}\footnotetext{
keywords: countable atomic boolean algebra, superatomic, small index property}
\newcounter{number}

\begin{abstract} It is shown that all the countable superatomic boolean algebras of finite rank have the small index property.
\end{abstract}

\section {Introduction}

In \cite{Dixon} it was shown that the full symmetric group on an infinite set has the `small index property' SIP, meaning that any 
subgroup having index strictly less than $2^{\aleph_0}$ contains the pointwise stabilizer of a finite set. The corresponding result 
for the group of order-preserving permutations of the set of rational numbers was given in \cite{Truss1}. The small index property 
has received a great deal of attention in quite a wide variety of special cases. Its model-theoretic significance is that its truth 
tells us that the natural topological group associated with the structure (under the topology of pointwise convergence) can be recovered 
from the pure group, and from this one can deduce that the structure is interpretable in the (abstract) automorphism group \cite{Hodges}. 
Most of the structures for which SIP has been studied are $\aleph_0$-categorical, but this is not required in the definition, and 
some non-$\aleph_0$-categorical cases have been looked at, for instance 1-transitive linear orders in \cite{Chicot}, and trees 
and cycle-free partial orders in \cite{Truss2}.

In this paper we look at a class of countable structures which are not $\aleph_0$-categorical, namely the countable superatomic boolean 
algebras of finite rank, and show that their automorphism groups have the small index property. The definitions here 
are that the boolean algebra $\mathbb B$ is {\em atomic} if every element is the least upper bound of the atoms below it, and it is 
{\em superatomic} if every homomorphic image is atomic. This notion was explored for instance in \cite{Mostowski} and \cite{Day}, 
and various equivalent conditions were given for superatomicity. Here we just look at the countable case (since for the most part, 
that is the context for considering SIP), where one can give an explicit description. The easiest one to work with is via topologies 
on countable ordinals. In fact the family of clopen subsets of any ordinal $\alpha$, which for technical reasons it is easiest to 
assume is a successor, forms a superatomic boolean algebra. Since for any successor ordinals $\alpha$ and $\beta$, $\alpha + \beta$ 
is homeomorphic to $\beta + \alpha$ under the order topology, by writing in Cantor normal form and reducing, the least successor 
ordinal in any homeomorphism class has the form $\omega^\alpha \cdot a + 1$ for some finite $a$, and in the terminology of \cite{Day}, 
it has cardinal sequence $(\aleph_0, \aleph_0, \ldots, a)$. 

To analyze countable superatomic boolean algebras $\mathbb B$, we consider the increasing sequence of ideals $(I_\beta: \beta \in On)$ 
given inductively as follows: $I_0 = \{0\}$; assuming $I_\beta$ has been defined, let $I_{\beta + 1} \supseteq I_\beta$ be such 
that $I_{\beta + 1}/I_\beta$ is the ideal of ${\mathbb B}/I_\beta$ generated by its set of atoms; for limits $\lambda$, 
$I_\lambda = \bigcup_{\beta < \lambda}I_\beta$. There is a least $\alpha = \alpha({\mathbb B})$ such that 
$I_\alpha = I_{\alpha + 1}$, and (by superatomicity), this can only happen when $I_\alpha = {\mathbb B}$. One can check that 
$\alpha({\mathbb B})$ is a successor, and $\alpha - 1$ is called the {\em rank} of $\mathbb B$. Clearly ${\mathbb B}/I_{\alpha -1}$ 
is finite, and ${\mathbb B}$ is determined uniquely up to isomorphism by its rank and the number $a$ of atoms of ${\mathbb B}/I_{\alpha - 1}$ 
(called its `degree'). From the point of view of verifying the small index property, we can assume that 
$|{\mathbb B}/I_{\alpha - 1}| = 2$ (so there is just one atom), since the general case can easily be derived from this (see Corollary \ref{3.2}). 
The ideal of a boolean algebra $\mathbb B$ generated by its atoms is called the {\em Fr\'echet ideal}, written $Fr({\mathbb B})$. Thus 
in the above inductive step, $I_{\beta + 1}/I_\beta = Fr({\mathbb B}/I_\beta)$.

Let us start then with an ordinal of the form $\omega^\alpha$, where this is the ordinal power, and we consider the successor 
ordinal $\omega^\alpha + 1$ under the ordinal topology. For technical reasons, we work with the interval $X = [1, \omega^\alpha]$ 
(which is order-isomorphic to $\omega^\alpha + 1$). If we let $\mathbb B$ be the boolean algebra of clopen subsets of $X$, 
then $\mathbb B$ is superatomic of rank $\alpha$ and degree 1. In fact one can verify that $I_\beta$ is the family of ordinals in 
this range which are not divisible by $\omega^{\beta + 1}$. For the proof that any countable superatomic boolean algebra is 
isomorphic to the algebra of clopen subsets of $\omega^\alpha \cdot a + 1$ for some countable ordinal $\alpha$ and finite $a$, 
see \cite{Day}. For simplicity we assume that $a = 1$, and deduce the result for general $a$ in the final section.

\section {The induction step}

We are aiming to prove that the group of homeomorphisms to itself $G = G_\delta$ of $[1, \delta]$ satisfies the small index property 
for all countable $\delta$ of the form $\omega^\alpha$. At present we have only succeeded in doing this for $\alpha$ finite.
This goes by induction on $\alpha$. Since however the induction step goes through for successors for arbitrary $\alpha$, the assumption 
that $\alpha$ is finite is only made in the final section. So we let $\delta = \omega^{\alpha + 1}$, where $\alpha$ is a countable 
ordinal, and this is the ordinal power. Let $X = X_\delta = [1, \delta]$, under the order topology, and let $\mathbb B$ be the 
boolean algebra of its clopen subsets. Then the group $G = G_\delta$ of homeomorphisms of $X$ to itself is the same as the 
automorphism group of $\mathbb B$. 

Let us write $A_n = [\omega^\alpha (n-1) + 1, \omega^\alpha \cdot n]$, for $n \in [1, \omega)$. Thus the sets $A_n$ form a 
partition of $X \setminus \{\delta\}$ into clopen subsets. Throughout this section we write $K$ for the subgroup of $G$ which fixes 
each $A_n$ setwise. We notice that each $A_n$ is homeomorphic to $[1, \omega^\alpha]$, so $K$ is isomorphic to the unrestricted 
direct product of the automorphism groups $G_n$ of $\mathbb B$ restricted to $A_n$, each of which is isomorphic to $G_{\omega^\alpha}$. 
For each $n$ we fix a bijection $\varphi_n$ from $A_1$ to $A_n$, which may be explicitly given by $\varphi_n(\beta) = \omega^{\alpha}(n-1) + \beta$
for $\beta \in [1, \omega^\alpha]$.  Then $\varphi_n$ is a homeomorphism, and so it induces an isomorphism from $G_1$ to $G_n$. 
Also $\varphi_{ij} = \varphi_j \varphi_i^{-1}$ is a homeomorphism from $A_i$ to $A_j$, and clearly $\varphi_{ij}^{-1} = \varphi_{ji}$.

\begin{lemma} If $N$ is a normal subgroup of $G$ of index $< 2^{\aleph_0}$, then $K \le N$. \label{2.1}  \end{lemma}

\noindent{\bf Proof}: We write $G_{\omega^\alpha} = H$, and identify $K$ with $H^\omega$ when desired. Let $[1, \omega)$ be written as 
the disjoint union $\bigcup_{n \in [1, \omega)}Z_n$ of infinite sets $Z_n$, and let $\theta_n$ be a bijection from $[1, \omega)$ to $Z_n$. 
Thus $\theta_{mn} = \theta_n\theta^{-1}_m$ is a bijection from $Z_m$ to $Z_n$ (and its inverse is $\theta_{nm}$). Let $h$ be 
an arbitrary member of $K$. Then we may `copy' $h$ to each $\bigcup_{i \in Z_n}A_i$ by letting 
$h_n = \bigcup_{i \in Z_n}\varphi_{i \, \theta_n^{-1}i}^{-1}h\varphi_{i \, \theta_n^{-1}i}$, and $h_n$ fixes all other points. 
Unravelling this, if $i \in Z_n$, then $\varphi_{i \, \theta_n^{-1}i}$ maps $A_i$ to $A_{\theta^{-1}_ni}$, where $h$ now acts, and 
then $\varphi_{i \, \theta_n^{-1}i}^{-1}$ maps this back to $A_i$. Note that the actions of all the $h_n$ are disjoint, since 
the support of $h_n$ is contained in $\bigcup_{i \in Z_n}A_i$.

Let $\cal M$ be a family of  $2^{\aleph_0}$ pairwise almost disjoint infinite subsets of $[1, \omega)$ such that the union of any 
two is coinfinite. For $M \in {\mathcal M}$, let $h_M$ agree with $h_n$ on $\bigcup_{i \in Z_n}A_i$ for all $n \in M$, and fix
all other members of $X$. Then $\{h_M: M \in {\mathcal M}\}$ is a family of $2^{\aleph_0}$ members of $G$, so two of them 
must lie in the same left coset of $N$. Hence $h_{M_1}^{-1}h_{M_2} \in N$ for some $M_1 \neq M_2$ in $\mathcal M$. Now if 
$n \in M_1 \cap M_2$, then $h_{M_1}$ and $h_{M_2}$ both agree with $h_n$ on $\bigcup\{A_i: i \in Z_n\}$, and so 
$h_{M_1}^{-1}h_{M_2}$ is the identity on $\bigcup\{A_i: i \in Z_n, n \in M_1 \cap M_2\}$, and it agrees with $h_n^{-1}$ 
on $\bigcup\{A_i: i \in Z_n\}$ for $n \in M_1 \setminus M_2$ and with $h_n$ on $\bigcup\{A_i: i \in Z_n\}$ for 
$n \in M_2 \setminus M_1$. We may therefore alternatively write $h_{M_1}^{-1}h_{M_2}$ as $h_{I_1}^{-1}h_{I_2}$ where 
$I_1 = M_1 \setminus M_2$ and $I_2 = M_2 \setminus M_1$ are disjoint and infinite, and with $I_1 \cup I_2$ coinfinite.

Now let $J_1$ and $J_2$ be any disjoint infinite subsets of $[1, \omega)$ with coinfinite union. Then there is a permutation $\psi$ 
of $[1, \omega)$ which takes $I_1$ to $J_1$ and $I_2$ to $J_2$. Let $k \in G$ be given by 
$k(x) = \varphi_{i \, \theta_{m \, \psi(m)}i}(x)$ if $x \in A_i$, $i \in Z_m$. Thus for each $m$, $k$ maps 
$\bigcup\{A_i: i \in Z_m\}$ bijectively to $\bigcup\{A_j: j \in Z_{\psi(m)}\}$. We want to show that 
$h_{J_1}^{-1}h_{J_2} = k h_{I_1}^{-1}h_{I_2} k^{-1}$. Take any $x \in X \setminus \{\delta\}$. Let $n \in \omega$ and $j \in Z_n$ be such that 
$x \in A_j$. Let $n = \psi(m)$, $\theta_{n m}j = i$, and $\varphi_{ij}(y) = x$. Then $i \in Z_m$, so $k$ maps $A_i$ to 
$\varphi_{i \, \theta_{m n} i}A_i = \varphi_{ij}A_i = A_j$, by $\varphi_{ij}$. Note also that from $\theta_{n m}j = i$ it follows 
that $\theta^{-1}_mi = \theta^{-1}_nj$. Suppose first that $m \in I_2$. Then $n \in J_2$. We calculate that 
\begin{eqnarray*}
kh_{I_1}^{-1}h_{I_2}k^{-1}(x) & = & \varphi_{ij}h_{I_1}^{-1}h_{I_2}\varphi_{ij}^{-1}(x) = \varphi_{ij}h_{I_1}^{-1}h_{I_2}(y) = \varphi_{ij}h_m(y)  \\
                              & = & \varphi_{ij}\varphi^{-1}_{i \, \theta_m^{-1} i}h \varphi_{i \, \theta^{-1}_m i} (y) = \varphi_{ij} \varphi_{\theta^{-1}_m i \, i}h \varphi_{i \, \theta_m^{-1}i}(y) \\
                              & = & \varphi_{\theta_m^{-1}i \, j}h \varphi_{i \, \theta_m^{-1} \, i} (y) = \varphi_{\theta^{-1}_n j \, j}h \varphi_{i \, \theta_n^{-1}j}(y) \\
                              & = & \varphi_{\theta_n^{-1}j \, j} h \varphi_{i \, \theta_n^{-1}j} \varphi_{ij}^{-1}(x) = \varphi_{\theta_n^{-1}j \, j}h\varphi_{i \, \theta_n^{-1} j} \varphi_{ji}(x) \\
                              & = & \varphi_{\theta_n^{-1} j \, j}h \varphi_{j \, \theta_n^{-1}j}(x) =
                               \varphi^{-1}_{j \, \theta_n^{-1} j}h \varphi_{j \, \theta_n^{-1} j}(x) = h_n(x) = h_{J_1}^{-1}h_{J_2}(x). 
\end{eqnarray*}
If $m \in I_1$, essentially the same calculation applies, with all elements involving $h$ replaced by their inverses. If 
$m \not \in I_1 \cup I_2$, then the calculation is as follows:
$$kh_{I_1}^{-1}h_{I_2}k^{-1}(x)  =  \varphi_{ij}h_{I_1}^{-1}h_{I_2}\varphi_{ij}^{-1}(x) = \varphi_{ij}h_{I_1}^{-1}h_{I_2}(y) = \varphi_{ij}(y)  =  x = h_{J_1}^{-1}h_{J_2}(x). $$
Since $N$ is normal, it follows that $h_{J_1}^{-1}h_{J_2}$ also lies in $N$. Take any such $J_1$ and $J_2$, and choose $J_1'$, 
$J_2'$ such that $J_2' \subseteq J_1$, $|J_1 \setminus J_2'| = 1$, and $J_1' = J_2$. By the above argument, there are members 
of $N$ which are equal to $h_{J_1}^{-1}h_{J_2}$ and $h_{J_1'}^{-1}h_{J_2'}$, and multiplying these, we get a member of $N$ which 
is equal to $h_{\{i\}}^{-1}$ for some $i \in \omega$. By applying further conjugacies as necessary, every $h_{\{j\}}$ lies in $N$.

The argument so far started with any $h \in K$, and found an infinite coinfinite set $Z$ such that any `copy' of $h$ with support 
contained in $Z$ lies in $N$. If instead we start with $h'$ having support contained in an infinite coinfinite set $Z_1$, then we can 
find a partition of $[1, \omega)$ into $Z_n$ for $n \in [1, \omega)$, and $h \in K$ such that $h_1 = h'$, which as we have just shown, lies in 
$N$. In other words, any member of $K$ whose support is contained in a coinfinite set lies in $N$. Now cut $\omega$ into two infinite 
pieces, and write $h = h'h''$ where these are the restrictions of $h$ to two pieces of 
the given kind. By the above argument, each of these restrictions lies in $N$, and hence so does $h$.     $\Box$

\vspace{.1in}

There is a rather larger subgroup of $G$ that we now need to consider, written $K^*$. This is defined to be the pointwise
stabilizer of $\{\omega^\alpha \cdot n: n \in [1, \omega]\}$. This clearly contains $K$, since $\omega^\alpha \cdot n$ is the unique
point of $A_n$ having Cantor--Bendixson rank $\alpha$. Again using the fact that all the $A_n$ are homeomorphic, we can permute them 
arbitrarily by members of $G$, and this gives a natural homomorphism $\pi$ of the setwise stabilizer of $\{A_n: 1 \le n < \omega\}$ 
onto ${\rm Sym}([1, \omega)) \;\; (\cong {\rm Sym}(\omega))$. If we extend this homomorphism to the whole of $G$ by letting 
$\pi(g)(m) = n$ for $m, n \in [1, \omega)$ whenever $g(\omega^\alpha \cdot m) = \omega^\alpha \cdot n$, then again
$\pi(g) \in {\rm Sym}[1, \omega)$, and we can identify $K^*$ as the kernel of $\pi$. 

Note that by continuity, whenever $g(\omega^\alpha \cdot m) = \omega^\alpha \cdot n$, there is some $x < \omega^\alpha \cdot m$ 
in $A_m$ such that $(x, \omega^\alpha \cdot m]$ is mapped into $A_n$. Thus the action of $G$ on $X$ is approximated by that of the setwise
stabilizer of $\{A_n: n \in [1, \omega)\}$. To handle this, we say that a subset of $A_n$ whose complement is bounded is 
{\em cofinal} (noting that this is not the usual sense of `cofinal'). A key part of our argument concerns how these bounded pieces 
are permuted by members of $G$.

\begin{lemma} If $N$ is a normal subgroup of $G$ of index $< 2^{\aleph_0}$, then for any $\sigma \in {\rm Sym} [1, \omega)$ there is 
$g \in N$ such that $\pi(g) = \sigma$. \label{2.2}  \end{lemma}

\noindent{\bf Proof}: Using the above notation we observe that $NK^*/K^*$ is a normal subgroup of $G/K^*$ of index $< 2^{\aleph_0}$, and 
since Sym $\,\omega$ has no proper normal subgroups of index $< 2^{\aleph_0}$, it follows that $NK^* = G$, and this establishes precisely 
what is wanted. For given $\sigma \in {\rm Sym} [1, \omega)$, let $g \in G$ be such that $\pi(g) = \sigma$ (for instance, $g$ may be 
taken in $K$). Since $NK^* = G$, we may write $g = hk$ where $h \in N$, $k \in K^*$, and it follows that $\pi(h) = \sigma$.  $\Box$

\vspace{.1in}

Let $\mathcal P$ be the family of ordered pairs of bounded clopen subsets of $X_{\omega^\alpha}$. If $P$ and $Q$ are topological spaces, 
we denote by $P \sqcup Q$ the space obtained by taking their disjoint union, with the topology under which each of $P$ and $Q$ is clopen, 
and the subspace topology on each of them is their original topology. If $(P_1, Q_1)$, $(P_2, Q_2) \in {\mathcal P}$, we write 
$(P_1 \sqcup P_2, Q_1 \sqcup Q_2)$ as $(P_1, Q_1) + (P_2, Q_2)$. There is an associated equivalence relation $\sim$ on $\mathcal P$ given 
by $(P_1, Q_1) \sim (P_2, Q_2)$ if $P_1 \sqcup Q_2$ is homeomorphic to $P_2 \sqcup Q_1$. Now for $g \in G$ we let 
$P_i = P_i(g) = A_{\pi(g)(i)} \setminus gA_i$, and $Q_i = Q_i(g) = A_i \setminus g^{-1}A_{\pi(g)(i)}$, so that 
$(P_i, Q_i) \in {\mathcal P}$ if we identify $X_{\omega^\alpha} = [1, \omega^\alpha]$ with $A_i$, since by definition of $\pi$, 
$g(\omega^\alpha \cdot i) = \omega^\alpha \cdot \pi(g)(i)$, and using continuity of $g$. In the first case one considers, which 
is $\alpha = 1$, $P_i$, $Q_i$ are finite, and up to $\sim$-equivalence we can just replace $(P_i(g), Q_i(g))$ by $|P_i(g)| - |Q_i(g)|$. In 
the general case, these sets are not necessarily finite, though they are bounded and clopen, and so in a sense they are encoded by a 
finite amount of information, being compact sets. If $(P, Q) \in {\mathcal P}$, we let $-(P, Q) = (Q, P)$, and this definition allows us to
subtract members of $\mathcal P$ (as well as add them).

\begin{lemma} If $g \in G$ and $B$ is a cofinal subset of $A_i$ such that $gB \subseteq A_{\pi(g)(i)}$, then 
$(P_i(g), Q_i(g)) \sim (A_{\pi(g)(i)} \setminus gB, A_i \setminus B)$. \label{2.3}  \end{lemma}

\noindent{\bf Proof}: Let us write $j$ for $\pi(g)(i)$. Then
\begin{eqnarray*}
(P_i(g), Q_i(g)) & = & (A_j \setminus gA_i, A_i \setminus g^{-1}A_j)  \\
                 & = & (A_j \setminus (A_j \cap gA_i), A_i \setminus (A_i \cap g^{-1}A_j))  \\  
                 & \sim & ((A_j \setminus (A_j \cap gA_i)) \cup ((A_j \cap gA_i) \setminus gB), (A_i \setminus (A_i \cap g^{-1}A_j)) \cup ((A_i \cap g^{-1}A_j) \setminus B))   \\
                 &   & \hspace{.5in} \mbox{ since } B \subseteq A_i \cap g^{-1}A_j \subseteq A_i, gB \subseteq A_j \cap gA_i \subseteq A_j, \\
                 &   & \hspace{.5in} \mbox{ and $g$ is a homeomorphism from } (A_i \cap g^{-1}A_j) \setminus B \mbox{ to } (A_j \cap gA_i) \setminus gB  \\
                 & = & (A_j \setminus gB, A_i \setminus B), \mbox{ as desired.} \hspace{.2in} \Box
\end{eqnarray*}

\begin{lemma} If $g, h \in G$ then for each $i$, $(P_i(hg), Q_i(hg)) \sim (P_i(g), Q_i(g)) + (P_{\pi(g)(i)}(h), Q_{\pi(g)(i)}(h))$. \label{2.4}  \end{lemma}

\noindent{\bf Proof}: Let us write $j = \pi(g)(i)$ and $k = \pi(h)(j)$ and by continuity of $g$ and $h$ pick a cofinal subset $B$ of $A_i$
such that $gB \subseteq A_j$ and $hgB \subseteq A_k$. Then by Lemma \ref{2.3},
\begin{eqnarray*}
(P_i(hg), Q_i(hg)) & = & (A_k \setminus hgB, A_i \setminus B) \\
        & \sim & ((A_k \setminus hgB) \sqcup (A_j \setminus gB), (A_j \setminus gB) \sqcup (A_i \setminus B))   \\
        & = & ((A_k \setminus hgB), (A_j \setminus gB)) + ((A_j \setminus gB), (A_i \setminus B))   \\
        & = & (P_i(g), Q_i(g)) + (P_j(h), Q_j(h)).  \hspace{.2in}  \Box
\end{eqnarray*}

\begin{lemma} For any $g$ such that $\sigma = \pi(g)$ maps $[1, \omega)$ in a single infinite cycle, and
$f: [1, \omega) \rightarrow {\mathcal P}$, there is $h \in K^*$ such that for each $i$, $(P_i(h^{-1}gh), Q_i(h^{-1}gh)) \sim f(i)$. 
\label{2.5}  \end{lemma}

\noindent{\bf Proof}: We choose $(R_i, S_i)$ inductively thus:
\begin{eqnarray*}
(R_1, S_1)                             & = & (\emptyset, \emptyset)                  \\
(R_{\sigma^{j+1}1}, S_{\sigma^{j+1}1}) & = & (R_{\sigma^j1}, S_{\sigma^j1}) - (P_{\sigma^j1}(g), Q_{\sigma^j1}(g)) + f(\sigma^j1) \mbox{ if } j \ge 0\\  
(R_{\sigma^{j-1}1}, S_{\sigma^{j-1}1}) & = & (R_{\sigma^j1}, S_{\sigma^j1}) + (P_{\sigma^{j-1}1}(g), Q_{\sigma^{j-1}1}(g)) - f(\sigma^{j-1}1) \mbox{ if } j \le 0. 
\end{eqnarray*}
Since $\sigma$ has a single infinite cycle, this defines $(R_i, S_i)$ for all $i \in [1, \omega)$, and we note
that $(R_{\sigma i}, S_{\sigma i}) - (R_i, S_i) = -(P_i(g), Q_i(g)) + f(i)$ for all $i$. For instance, if $i = \sigma^j1$ where
$j < 0$, $(R_{\sigma^j1}, S_{\sigma^j1}) = (R_{\sigma^{j+1}1}, S_{\sigma^{j+1}1}) + (P_{\sigma^j1}(g), Q_{\sigma^j1}(g)) - f(\sigma^j1)$,
which gives $(R_{\sigma i}, S_{\sigma i}) - (R_i, S_i) = -(P_i(g), Q_i(g)) + f(i)$.

Now for each $i$ we shall choose bounded clopen subsets $B_i, C_i$ of $A_i$ such that $(B_i, C_i) \sim (R_i, S_i)$,
and show that (by careful choice of $B_i$ and $C_i$) there is $h \in K^*$ which maps $A_i \setminus B_i$ onto $A_i \setminus C_i$. The fact that 
$A_i \setminus B_i$ can be taken to $A_i \setminus C_i$ is immediate, since they are even order-isomorphic, and being clopen, we can act on
the complement $\bigcup_{i \in \omega}B_i$ as we please. If this has been done, we can see that $(P_i(h), Q_i(h)) \sim -(R_i, S_i)$, and 
similarly $(P_i(h^{-1}), Q_i(h^{-1})) \sim (R_i, S_i)$. For this note first that since $h \in K^*$, $\pi(h) = id$. Hence
\begin{eqnarray*}
(P_i(h), Q_i(h)) & = & (A_i \setminus hA_i, A_i \setminus h^{-1}A_i) \\
        & = & (A_i \setminus (h(A_i \setminus B_i) \cup hB_i), A_i \setminus (h^{-1}(A_i \setminus C_i) \cup h^{-1}C_i))   \\
        & = & (C_i \setminus hB_i, B_i \setminus h^{-1}C_i)   \\
        & \sim & (C_i, B_i) = -(R_i, S_i). 
\end{eqnarray*}
To justify the final step, we have to see that $(C_i \setminus hB_i) \sqcup B_i$ is homeomorphic to $(B_i \setminus h^{-1}C_i) \sqcup C_i$. Now, 
$(C_i \setminus hB_i) \sqcup B_i$ is homeomorphic to $(h^{-1}C_i \setminus B_i) \sqcup B_i = B_i \cup h^{-1}C_i = 
(B_i \setminus h^{-1}C_i) \cup h^{-1}C_i$, which is homeomorphic to $(B_i \setminus h^{-1}C_i) \sqcup C_i$ as desired.

Fix $i$ and let $\sigma i = j$. Then by Lemma \ref{2.4} we find that $(P_i(h^{-1}gh), Q_i(h^{-1}gh)) \sim (P_i(h), Q_i(h)) + (P_i(g), Q_i(g)) + (P_j(h^{-1}), Q_j(h^{-1})) 
\sim -(R_i, S_i) + (P_i(g), Q_i(g)) + (R_{\sigma i}, S_{\sigma i}) \sim f(i)$, as desired. 

To choose $B_i$ and $C_i$ we note that in the case $\alpha = 1$, any non-empty bounded clopen sets will serve, and this is because 
they will be finite, and so $\bigcup_{i \in [1, \omega)}B_i$ and $\bigcup_{i \in [1, \omega)}C_i$ are (countably) infinite sets with the discrete topology, so any 
bijection from the first to the second works. In the general case we have to take the topology into account. For this, let
$B_i' = B_{i1} \sqcup B_{i2} \sqcup \ldots \sqcup B_{ii} \sqcup C_{i1} \sqcup C_{i2} \sqcup \ldots \sqcup C_{i \, i-1}$
and $C_i' = B'_{i1} \sqcup B'_{i2} \sqcup \ldots \sqcup B'_{i \, i-1} \sqcup C'_{i1} \sqcup C'_{i2} \sqcup \ldots \sqcup C'_{ii}$,
where $B_{ij}$, $B'_{ij}$ are homeomorphic to $B_j$, and $C_{ij}$, $C'_{ij}$ to $C_j$. Then $(B_i', C_i') \sim (B_{ii}, C_{ii})
\sim (B_i, C_i) \sim (R_i, S_i)$, and $B_i'$ and $C_i'$ are still homeomorphic to bounded clopen subsets of $X_\alpha$, so will do equally well 
for our definition. We can now map the disjoint union of all the $B_i'$s to that of the $C_i'$s by taking $B_{ij}$ to $B'_{i+1 \, j}$ and 
$C_{i+1 \, j}$ to $C'_{i j}$ where each individual map is a homeomorphism, and one verifies that it takes $\bigsqcup_{i \in [1, \omega)}B_i'$ 
1--1 onto $\bigsqcup_{i \in [1, \omega)}C_i'$, so is a suitable choice for our extension of $h$.    $\Box$

\begin{lemma} If there is no proper normal subgroup of $G_{\omega^\alpha}$ of index $< 2^{\aleph_0}$, then the same is true for
$G = G_{\omega^{\alpha+1}}$. \label{2.6}  \end{lemma}

\noindent{\bf Proof}: Suppose that $N \lhd G$ and $|G: N| < 2^{\aleph_0}$, and we shall show that $N = G$. Let $g \in G$ 
be arbitrary. If $\sigma$ permutes $\omega$ in a single cycle, then by Lemma \ref{2.2} there is $h \in N$ such that 
$\pi(h) = (\pi(g))^{-1}\sigma$. Thus $\pi(gh) = \pi(g)\pi(h) = \sigma$ is a single infinite cycle, and applying Lemma \ref{2.2} again, 
there is $k \in N$ such that $\pi(k) = \pi(gh)$. By Lemma \ref{2.5} we may replace $k$ by a conjugate and suppose that 
$(P_i(k), Q_i(k)) \sim (P_i(gh), Q_i(gh))$ for each $i$ (noting that the conjugating element has trivial $\pi$ value). 
By Lemma \ref{2.4}, 
\begin{eqnarray*}
(P_i(k^{-1}gh), Q_i(k^{-1}gh)) & \sim &  (P_i(gh), Q_i(gh)) + (P_{\sigma i}(k^{-1}), Q_{\sigma i}(k^{-1}))                    \\
              & \sim &  (P_i(k), Q_i(k)) + (P_{\sigma i}(k^{-1}), Q_{\sigma i}(k^{-1}))     \\  
              & \sim &  (P_i(k^{-1} \cdot k), Q_i(k^{-1} \cdot k)) \sim (\emptyset, \emptyset).
\end{eqnarray*}
Thus $k^{-1}gh$ is a member of $K^*$ such that $(P_i(k^{-1}gh), Q_i(k^{-1}gh)) \sim (\emptyset, \emptyset)$ for all $i$. Let 
$B_i = \{a \in A_i: k^{-1}gha \not \in A_i\}$ and $C_i = \{a \in A_i: (k^{-1}gh)^{-1}a \not \in A_i\}$. Thus $B_i$ and
$C_i$ are homeomorphic bounded clopen subsets of $A_i$. Furthermore, 
$A_i \setminus C_i = \{a \in A_i: (k^{-1}gh)^{-1}a \in A_i\} = (k^{-1}gh)A_i \cap A_i = (k^{-1}gh)(A_i \setminus B_i)$.
Let $l \in G$ agree with $k^{-1}gh$ on $A_i \setminus B_i$ and map $B_i$ to $C_i$. 
Thus $l^{-1}k^{-1}gh$ fixes $\bigcup_{i \in \omega}(A_i \setminus B_i)$ pointwise, and
permutes points of $\bigcup_{i \in \omega}B_i$. Also by Lemma \ref{2.1}, $l \in K \le N$.

Since each $B_i$ is a bounded subset of $A_i$, $\bigcup_{i \in \omega}B_i$ has order-type $\le \omega^\alpha$, so there are 
bounded clopen subsets $D_i$ of $A_i$ such that $B_i \subseteq D_i$ and $\bigcup_{i \in \omega}D_i \cong \omega^\alpha$. Now
$Y = \bigcup_{i \in \omega}D_i$ is open, and has open complement, so is clopen. Let $L$ be the subgroup of $G$ comprising its 
elements whose support is contained in $Y$. Since $Y$ is clopen, any homeomorphism of $Y$ to itself extends to a
homeomorphism of $X$ by fixing $X \setminus Y$ pointwise, and therefore $L$ is a subgroup of $G$ which is isomorphic to 
$G_{\omega^\alpha}$. Since $N \cap L$ is a normal subgroup of $L$ of index $< 2^{\aleph_0}$, by assumption it follows that 
$L \le N$. Since $l^{-1}k^{-1}gh$ fixes $\bigcup_{i \in \omega}(A_i \setminus B_i)$ pointwise it also fixes
$\bigcup_{i \in \omega}(A_i \setminus D_i)$ pointwise, so lies in $L$, and hence also in $N$. Hence 
$g = kl(l^{-1}k^{-1}gh)h^{-1} \in N$ as required. \hspace{.2in} $\Box$

\begin{corollary}  For finite $\alpha$, $G_{\omega^\alpha}$ has no proper normal subgroup of index $< 2^{\aleph_0}$.
  \label{2.7}   \end{corollary}
  
\noindent{\bf Proof}: This follows by induction from the lemma. In the basis case $\alpha = 1$,  
$G_{\omega^\alpha} \cong {\rm Sym}(\omega)$, whose non-trivial normal subgroups are explicitly known, and are the 
alternating group, and the group of elements of finite support, both of which have index $2^{\aleph_0}$.
$\Box$

\section{The main result}

We now use the results of section 2 to complete the main proof.

\begin{theorem} The boolean algebra $\mathbb B$ of clopen subsets of $[1, \omega^\alpha]$ for finite $\alpha$ has the small 
index property. This is, any subgroup $H$ of its automorphism group $G$ of index $< 2^{\aleph_0}$ contains
the pointwise stablilizer $G_A$ of some finite $A \subseteq {\mathbb B}$. \label{3.1} \end{theorem}

\noindent{\bf Proof}: We use induction. The basis case $\alpha = 1$ follows from the small index property for
${\rm Sym}(\omega)$ \cite{Dixon}. So now assume the result for $\alpha$, and we prove it for $\alpha + 1$. So 
$G = G_{\omega^{\alpha + 1}}$, and as usual we let $A_n = [\omega^\alpha (n-1) + 1, \omega^\alpha \cdot n$. 
Let $K$ be the subgroup of $G$ fixing $\{A_n: n \in \omega\}$ setwise (not the same $K$ as before). We follow the proof 
of Lemma \ref{2.1}. As there, let $\pi: G \to {\rm Sym}[1, \omega)$ where $\pi(g)(m) = n$ if 
$g(\omega^\alpha \cdot m) = \omega^\alpha \cdot n$. Then $K^*$ given above is the kernel of $\pi$. 

Since $|G: H| < 2^{\aleph_0}$, we deduce that $|K: K \cap H| < 2^{\aleph_0}$, and also $|\pi K: \pi(K \cap H)| < 2^{\aleph_0}$. 
Now $\pi K$ is naturally isomorphic to ${\rm Sym}[1, \omega)$, which by \cite{Dixon} has the small index property.
Hence there is a finite $I \subseteq [1, \omega)$ such that any member $\sigma$ of ${\rm Sym}[1, \omega)$ fixing $I$ pointwise
lies in $\pi(K \cap H)$, so has the form $\pi(g)$ for some $g \in K \cap H$.

Let $G_1$ and $G_2$ be the restrictions of $G$ to $\bigcup_{i \in I}A_i$ and $\bigcup_{i \not \in I}A_i$ respectively. That is, 
$G_1$, $G_2$ consist of all members of $G$ whose support is contained in $\bigcup_{i \in I}A_i, \bigcup_{i \not \in I}A_i$ 
respectively. For each $i$, let $K_i = \{g \in G: supp(g) \subseteq A_i\}$. Then $K_i \cong G_{\omega^\alpha}$ and 
$|K_i: K_i \cap H| < 2^{\aleph_0}$, so by the induction hypothesis, there is a finite family $B_i$ of clopen subsets
of $A_i$ such that $(K_i)_{B_i} \le H$. Hence $H \ge \prod_{i \in I}(K_i)_{B_i} = (G_1)_{\bigcup_{i \in I}B_i \, \cup \{A_i: i \in I\}}$. 
We let $Y = \bigcup_{i \in I} B_i \cup \{A_i: i \in I\}$ and show that $H \ge G_Y$. Since $G_Y = (G_1)_Y \times G_2$, it 
remains to see that $G_2 \le H$. This follows by standard arguments as in \cite{Dixon} which we briefly sketch.

First write $[1, \omega) \setminus I$ as the disjoint union of infinite sets $M_i$ for $i \in [1, \omega)$. By considering the 
projections to $L_i =$ the restriction of $G$ to $\bigcup_{j \in M_i}A_j$ we find as in \cite{Dixon} that $H$ projects 
onto some $L_i$. It follows that $H \cap L_i$ is a normal subgroup of $L_i$ of index $< 2^{\aleph_0}$. Since $L_i \cong G$, 
by Corollary \ref{2.7}, $L_i \le H$. By conjugating by elements of $H$ arbitrarily permuting the members of $\omega \setminus I$ 
we deduce that $H$ contains the restriction of $G$ to $\bigcup_{j \in J}A_j$ for an arbitrary infinite coinfinite subset $J$ 
of $\omega \setminus I$. Since $G_2$ is generated by such restrictions, it follows that $H \ge G_2$.        $\Box$

\begin{corollary} Any superatomic boolean algebra $\mathbb B$ of finite rank has the small index property. \label{3.2}  \end{corollary}

\noindent{\bf Proof}: Let $\mathbb B$ have rank $\alpha$ and degree $a$. In the first case, if $a = 1$, the result is read off 
from the theorem. For the second part, let $\mathbb B$ be isomorphic to the family of clopen subsets of $[1, \omega^\alpha \cdot a]$, 
and for $1 \le i \le a$ let $A_i = [\omega^\alpha (i-1) + 1, \omega^\alpha \cdot i]$. Note that each $A_i$ is clopen, so lies in
$\mathbb B$. Furthermore, $G_{\{A_1, \ldots, A_a\}}$ is the direct product of its restrictions $G_i$ to the individual $A_i$s,
each of which is isomorphic to $G_{\omega^\alpha}$. Now $|G_i: G_i \cap H| < 2^{\aleph_0}$, so by the theorem, 
$G_i \cap H \ge (G_i)_{B_i}$ for some finite set $B_i$ of clopen subsets of $A_i$. Hence 
$H \ge G_{B_1 \cup \ldots \cup B_a \cup \{A_i:1 \le i \le a\}}$, concluding the proof.  \hspace{.2in}   $\Box$

\vspace{.2in}

\noindent{\bf Concluding remarks}

The main problem left open in what we have done is to prove the small index property for the homeomorphism group of $[1, \omega^\alpha]$
for an arbitrary countable ordinal $\alpha$, which by the argument of Corollary \ref{3.2} would give SIP for all countable 
superatomic boolean algebras. As we have seen it would suffice to be able to handle the limit step. The methods given
here do not seem to apply to this situation. We remark that there is a related piece of work in the final chapter of \cite{Hilton}, in 
which Hilton actually classified {\em all} the normal subgroups of $G_{\omega^\alpha \cdot a}$ again for finite $\alpha$, of which there are 
$2^{2^{\aleph_0}}$ (provided $\alpha \ge 2$). For what we require here, all that was needed was to know that any such proper normal 
subgroup has index $2^{\aleph_0}$.

\end{document}